\documentclass[10 pt,a4paper]{article}
\usepackage{latexsym,amssymb,graphicx,amscd,amsmath,url}
\usepackage[french,english]{babel}

\title{Limits of the quantum $SO(3)$ representations for the one-holed torus}

\author{Ramanujan Santharoubane}

\date{}

\begin{document}

\maketitle

\begin{abstract} For $N  \geq 2$, we study a certain sequence $(\rho_p^{(c_p)})$ of N-dimensional representations of the mapping class group of the one-holed torus arising from $SO(3)$-TQFT, and show that the conjecture of
Andersen, Masbaum, and Ueno \cite{1} holds for these representations. This is done by proving that, in a certain basis and up to a rescaling, the matrices of these representations converge as $p$ tends to infinity. Moreover, the limits describe the action of $SL_{2}(\mathbb{Z})$ on the space of homogeneous polynomials of two variables of total degree $N-1$.
\end{abstract}

\section{Introduction}

Quantum topology was born by the physical interpretation of the Jones Polynomials made by E. Witten. An interesting problem in quantum topology is to study the asymptotics of quantum objects by linking them to classical objects. In this paper we focus on the quantum representations of the mapping class group of the one-holed torus arising from the Witten-Reshetikhin-Turaev $SO(3)$ Topological Quantum Field Theory (TQFT).

For any  odd $p$ and any $c \in \lbrace 0,...,\frac{p-3}{2} \rbrace$, the $SO(3)$ TQFT built in [2] associates to the one-holed torus a finite dimensional complex vector space $V_p(T_{c})$ of dimension $\frac{p-1}{2}-c$. Denoting $\Gamma_{1,1}$ the mapping class group of the one-holed torus, $V_p(T_c)$  carries a projective representation of $\Gamma_{1,1}$ which depends on a choice of a primitive 2p-th root of unity $A_p$. It is known (see \cite{4}) and easy to see that in the case of $\Gamma_{1,1}$, this projective represention lifts to a linear representation which we denote:
\[ \rho_{p}^{(c)} :  \Gamma_{1,1} \longrightarrow Aut(V_p(T_{c}))\]

On the other hand, $\Gamma_{1,1}$ maps onto $SL_2(\mathbb{Z})$. For $N \geq 2$, the later group acts naturally on $H_N$ : the space of homogeneous polynomials of two variables of total degree $N-1$. So we have a representation:
\[ h_N :  \Gamma_{1,1} \longrightarrow Aut(H_N)\]

Remark that if p is odd and $p \geq 2N+1$ we can set $ c_p = \frac{p-1}{2}-N$ so that $\dim (V_p(T_{c_p}))= N$. This creates a sequence of N-dimensional representations $\rho_{p}^{(c_p)}$ of $ \Gamma_{1,1}$. It turns out that those representations are closely related to $h_N$. Indeed, up to rescaling, the quantum representations can be viewed as deformations of $h_N$. Here is a precise statement of what we mean:

\paragraph{Main theorem :}Let $\mathbb{Q}(X)$ be the field of rational functions in an indeterminate $X$. Fix $N \geq 2$ an integer. There exists a representation $\rho : \Gamma_{1,1}   \longrightarrow GL_N(\mathbb{Q}(X))$ which does not depend on $p$ and a character $\chi_p : \Gamma_{1,1} \longrightarrow \mathbb{C}^{*}$ (which depends on the choice of root of unity $A_p$)  such that :
\begin{itemize}
\item All the matrices in $ \rho(\Gamma_{1,1}) $ can be evaluated at $X= A_p$ and $X=-1$, those evaluations are denoted respectively $\rho^{[A_p]}$ and $\rho^{[-1]}$ (which are representations into $GL_{N}(\mathbb{C})$)

\item $\chi_p  \otimes  \rho^{[A_p]}$ is isomorphic to $\rho_p^{(c_p)}$

\item $ \rho^{[-1]}$ is isomorphic to $h_N$
\end{itemize}

\paragraph{} Let $t_{y}$ and $t_{z}$ be the Dehn twists along the canonical meridian and longitude on the one-holed torus. We choose the map $\Gamma_{1,1} \rightarrow SL_{2}(\mathbb{Z})$ such that $t_y$ maps to $ \begin{pmatrix}
1& 1 \\
0&1
\end{pmatrix}$ and  $t_{z}$ maps to $\begin{pmatrix}
1& 0 \\
-1&1
\end{pmatrix}$. Since  $t_{y}$ and $t_{z}$ generate $\Gamma_{1,1}$, the main theorem is implied by the following one:

\paragraph{Theorem 1 :}
Let  $\lbrace Q_n^{'(c_p)}\rbrace_{0\leqslant n \leqslant N-1}$ be the basis of $V_p(T_c)$ defined in \cite{3}.  Let $T_p$ and $T_p^{*}$ be the matrices of $\mu_{c_p}^{-1}\rho_{p}^{(c_p)}(t_y)$ and $\mu_{c_p}^{-1}\rho_{p}^{(c_p)}(t_{z})$ in this basis, where $\mu_{c_p}= (-A_p)^{c_p(c_p+2)}$. Then there exists  $\widehat{T}(X) , \widehat{T}^{*}(X) \in GL_N(\mathbb{Q}(X))$ independent of $p$ which can be evaluated at $X= A_p$ and $X=-1$ such that:
\begin{itemize}
\item $T_p = \widehat{T}(A_p)$ and  $T^{*}_p = \widehat{T}^{*}(A_p)$
\item  The matrices $\widehat{T}(-1)$ and $\widehat{T}^{*}(-1)$ are the matrices of $h_N(t_{y})$ and $h_N(t_{z})$ in the basis
$$\alpha_{n}X^{N-n-1}Y^{n}\ \ , \ \ 0\leq n \leq N-1  \leqno(1) $$of $H_{N}$, where $\alpha_{n}=\dfrac{2^{n}}{n!(N-1-n)!}$.

\end{itemize}

\paragraph{Remark :}Concretely, the previous theorem implies that if $\phi \in \Gamma_{1,1}$ and $M_p$ denotes the matrix of $\rho_p^{(c_p)}(\phi)$ (in the basis of theorem 1), we have as $A_p \rightarrow -1$:

\[ (\chi_p(\phi))^{-1} M_{p} \rightarrow M \] where $M$ is the matrix of $h_N(\phi)$ in the basis $(1)$.

\paragraph{} We can also use this theorem to prove the following version of the AMU  conjecture (see [1]) in the case of the one-holed torus:

\paragraph{Theorem 2 :}For any fixed $N \geq 2$, if $\phi \in \Gamma_{1,1}$ is pseudo-Anosov then there exists $p_0(\phi)$ such that for any odd $p\geq p_0(\phi)$ the automorphism $\rho_{p}^{(c_p)}(\phi)$ has infinite order.

\paragraph{Remark :}The case $N=2$ was already known to G. Masbaum ; see [1, Remark 5.9] and [3, p.96].

\paragraph{\textbf{Acknowlegments.}} I would like to thank Gregor Masbaum who gave me this problem and who helped me to write a precise statement of my result.

\section{Review of $SO(3)$-TQFT}

We are going to recall the basic notions we need. We refer to \cite{3} for more details.
\subsection{Notations}
For $p$ an odd integer, let $A=A_p$ be a primitive 2p-th root of unity. Let $n$ be an integer, we put $\lbrace n \rbrace=(-A)^{n}-(-A)^{-n}$, $\lbrace n \rbrace^{+}=(-A)^{n}+(-A)^{-n}$. When n is positive let $\lbrace n \rbrace! = \lbrace 1 \rbrace...\lbrace n \rbrace$ with $\lbrace 0 \rbrace!=1$ and when n is negative $\lbrace n \rbrace! =0$. Also $\lbrace n \rbrace!! = \lbrace n \rbrace \lbrace n-2 \rbrace...$. We put $\mu_n= (-A)^{n(n+2)}$ and $\lambda_{n}= -\lbrace 2n+2 \rbrace$.

In what follows, $N \geq 2$ will be fixed. We set $c=d-N$ where $d=\frac{p-1}{2}$. We saw that $V_p(T_c)$ is N-dimensional and it has a basis $\lbrace L_{c,n} \rbrace_{0 \leq n \leq N-1}$ given by the colored graphs in the solid torus (see \cite{3}) which can be described pictorially by the following diagrams :

$ L_{c,n}  \  \ 
=
  \  \ 
\begin{minipage}{0.9in}\includegraphics[width=0.9in]{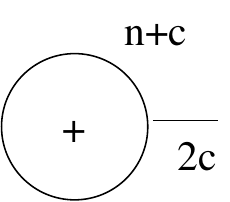} \end{minipage}
$

\paragraph{} $T_c$ can be viewed as a torus $\mathbb{T}^{2}$ equipped with a banded point $x$ with color $2c$. We can think of $T_c$ as the boundary of a tubular neighborhood of the graph.This tubular neighborhood is a solid torus, and the univalent vertex of the graph is "attached" to the banded point $x$.

\paragraph{} $\rho_p^{(c)}(t_{y})$ has a nice expression in this basis, indeed  for $ 0  \leq n \leq N-1 $
\[ \rho_p^{(c)}(t_{y})(L_{c,n}) = \mu_{c+n} L_{c,n} \]

We also denote by $((,))$ the Hopf pairing on $V_p(T_c)$. It is a symmetric non-degenerated bilinear form (see \cite{3}).

\paragraph{Remark :} In [3], $\rho_p^{(c)}(t_{y})$ is denoted by $t$ and   $\rho_p^{(c)}(t_{z})$ is denoted by $t^{*}$.

\subsection{Curve operators}
 For any multicurve (disjoint union of simple close curves) $\gamma$ on the one-holed torus $\mathbb{T}^{2}-x$ we can form the cobordism $C_{\gamma}$ as $\mathbb{T}^{2} \times I$ (where $I= [0,1]$) equipped with the banded link $ \gamma \times [\frac{1}{2},\frac{3}{4}]\cup x \times I$ where $x \times I$ has color $2c$. By the axioms of TQFT, $C_{\gamma}$ defines an operator $Z_p(\gamma) \in End(V_p(T_c))$. Let $y$ and $z$ be respectively the meridian and the longitude curves on the one-holed torus. We can see the action of $Z_p(y)$ and $Z_p(z)$ in the basis $\lbrace L_{c,n} \rbrace$:

 $ Z_p(y)(L_{c,n})  \  \ 
=
  \  \ 
\begin{minipage}{0.9in}\includegraphics[width=0.9in]{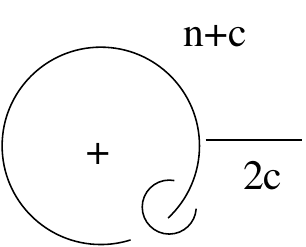} \end{minipage}
$

 $ Z_p(z)(L_{c,n})  \  \ 
=
  \  \ 
\begin{minipage}{0.9in}\includegraphics[width=0.9in]{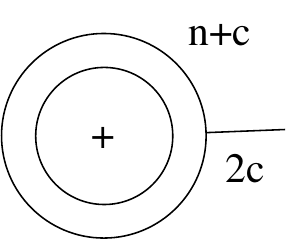} \end{minipage}
$

\paragraph{}Those diagrams can be evaluated using skein theory. We have also the following basic facts:

\paragraph{Proposition : }
Let $\phi \in \Gamma_{1,1}$ and $\gamma$ be a simple closed curve.
\begin{itemize}
\item $\rho_p^{(c)}(\phi) Z_p(\gamma) (\rho_p^{(c)}(\phi))^{-1} = Z_p(\phi(\gamma)) $
\item $Z_p(\gamma)$ is diagonisable and it has $\lbrace \lambda_c,...,\lambda_{c+N-1} \rbrace$ as eigenvalues.
\item $Z_p(z)$ and $Z_p(y)$ are transposed by the Hopf pairing.
\item $Z_p(t_{y}(z))= \dfrac{AZ_p(y)Z_p(z)- A^{-1}Z_p(z)Z_p(y)}{\lbrace 2 \rbrace}$
\end{itemize}

\paragraph{Remark :}The last property is obtained by applying the skein relation:

\[ 
\begin{minipage}{0.4in}\includegraphics[width=0.4in]{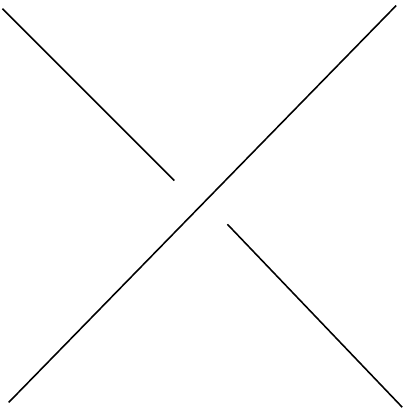}\end{minipage}
=  
 \ \ A \  \ 
\begin{minipage}{0.4in}\includegraphics[width=0.4in]{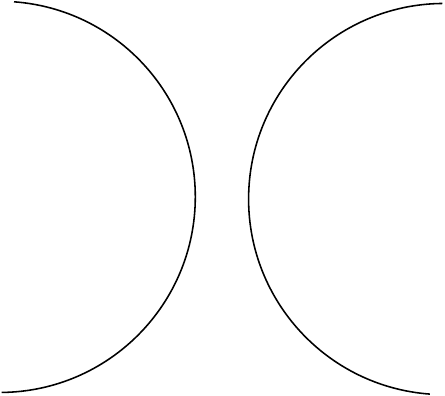} \end{minipage}
+
\ \  A^{-1} \  \ \begin{minipage}{0.4in} \includegraphics[width=0.4in]{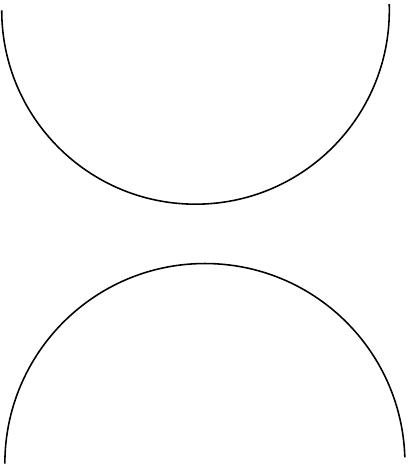}  \end{minipage}
   \]

\subsection{The basis $ \lbrace Q_n^{'(c)} \rbrace $}

We can now recall the definition of the basis used in theorem 1. Following \cite{3}, for $0 \leq n \leq N-1$ let:
$$  Q_n^{'(c)}= (\lbrace n \rbrace!)^{-1} \Biggl( \prod _{j=0}^{n-1}(Z_p(z)-\lambda_{c+j}Id) \Biggr)\big(L_{c,0}\big) $$

The interest of this basis is that it is orthogonal with respect to the Hopf pairing. For $n,m$ we have (see [3]):
\[ \dfrac{((Q_n^{'(c)},Q_n^{'(c)}))}{((Q_m^{'(c)},Q_m^{'(c)}))} = \dfrac{\lbrace m \rbrace ! \lbrace 2c+2n+1 \rbrace !! \lbrace 2c+n+1 \rbrace^{+}!}{\lbrace n \rbrace ! \lbrace 2c+2m+1 \rbrace !! \lbrace 2c+m+1 \rbrace^{+}!} \]

This quantity will be denoted by $R_{n,m}^{(c)}$. In the basis $ \lbrace Q_n^{'(c)} \rbrace$, let $T_p=(a_{m,n})$ be the matrix of $\mu_c^{-1}\rho_{p}^{(c)}(t_{y})$ and $T^{*}_p=(b_{n,m})$ be the matrix of $\mu_c^{-1}\rho_{p}^{(c)}(t_{z})$.
Let also $(y_{m,n})$ be the matrix of $Z_p(y)$, $(z_{m,n})$ be the matrix of $Z_p(z)$ and $(z^{'}_{m,n})$ be the matrix of $Z_p(t_y(z))$ . Since, with respect to the Hopf pairing,  $\rho_{p}^{(c)}(t_y)$ is the transpose of $\rho_{p}^{(c)}(t_{z})$ and $Z_p(y)$ is the transpose of $Z_p(z)$,we have :  

\[ a_{m,n} = R_{n,m}^{(c)} b_{n,m} \] and \[  y_{m,n} = R_{n,m}^{(c)} z_{n,m} \]

\paragraph{Remark : } In \cite{3} there are already explicit expressions for $(a_{m,n})$ and $(b_{n,m})$ but here we use new formulas which are more helpful for our purpose.

\section{The limit of the representations}

In this section we will prove Theorem 1 and Theorem 2.

\subsection{Proof of theorem 1}

 To make the proof of this theorem easier we need the following two lemmas:

\paragraph{Lemma 1.} 
For any $n$ and $m$, there exists $\widehat{R}_{n,m}(X) \in \mathbb{Q}(X)$ independent of $p$ which can evaluated at $X=A$ and $X=-1$ such that :

\begin{itemize}
\item $R_{n,m}=\widehat{R}_{n,m}(A) $

\item $\widehat{R}_{n,m}(-1) = \dfrac{(-4)^{n-m} m!(N-1-m)!}{n!(N-1-n)!}$

\end{itemize}

\subparagraph{Proof:} If $n=m$ the result is clear. By symmetries, we just have to prove it for $n \geq m+1$. 
 In this case, since $(-A)^p = 1$ and $2c=p-1-2N$ we have for any $x$ integer:
\[\begin{aligned}
\lbrace x+2c \rbrace   \ \ = \ \ \lbrace x-1-2N \rbrace &=- \lbrace2N+1-x \rbrace \\
 \lbrace x+2c \rbrace^{+}=  \lbrace x-1-2N \rbrace^{+} &= \ \  \lbrace2N+1-x \rbrace^{+}
\end{aligned}\] so\[\begin{aligned}
R_{n,m}&= \dfrac{\lbrace m \rbrace ! \lbrace 2c+2n+1 \rbrace !! \lbrace 2c+n+1 \rbrace^{+}!}{\lbrace n \rbrace ! \lbrace 2c+2m+1 \rbrace !! \lbrace 2c+m+1 \rbrace^{+}!} \\
&= \frac{(-1)^{n-m} \lbrace m \rbrace ! \lbrace 2N-2m-2 \rbrace !! \lbrace 2N-m-1 \rbrace ^{+}!}{\lbrace n \rbrace ! \lbrace 2N-2n-2 \rbrace !! \lbrace 2N-n-1 \rbrace ^{+}!} 
\end{aligned}\]Since for all $a$, $\lbrace a \rbrace$ is the evaluation of $(-X)^{n}-(-X)^{-n} \in \mathbb{Q}(X)$ at $X=A$ and $\lbrace a \rbrace^{+}$ is the evaluation of $(-X)^{n}+(-X)^{-n} \in \mathbb{Q}(X)$ at $X=A$ we see that there exists $\widehat{R}_{n,m}(X) \in \mathbb{Q}(X)$ (which clearly does not depend on $p$) such that $R_{n,m}=\widehat{R}_{n,m}(A)$. We also know that for all $a$ , when $A$ tends to $-1$ , $\frac{\lbrace a \rbrace}{\lbrace 1 \rbrace} \rightarrow a$ and $ \lbrace  a \rbrace^{+} \rightarrow 2$ so we deduce the expression of  $\widehat{R}_{n,m}(-1)$  $\square$

\paragraph{Lemma 2.} 

For $n \leq N-2$, let $M^{(n)}=(M_{m,l}^{(n)})$ be the matrix of $$\lbrace n+1 \rbrace^{-1}(Z_{p}(t_{y}(z))-\lambda_{c+n}Id_N)$$ in the basis $ \lbrace Q_n^{'(c)} \rbrace$.

 Then there exists a matrix $\widehat{M}^{(n)}(X)=(\widehat{M}^{(n)}_{m,l}(X)) \in GL_N(\mathbb{Q}(X)) $   independent of $p$ which can be evaluated at $X=A$ and $X=-1$ such that for all $m,l$:

\begin{itemize}

\item $\widehat{M}^{(n)}_{m,l}(A)= M_{m,l}^{(n)}$

\item $\widehat{M}^{(n)}_{m,l}(-1)=  \delta_{m-1,l}\dfrac{m}{n+1}+\delta_{m,l}\dfrac{2(N-2m-1)}{n+1} +\delta_{m+1,l}\dfrac{-4(N-m-1)}{n+1}$

\end{itemize}
where $\delta_{k,l}$ is the Kronecker symbol.

\subparagraph{Proof :} Since $Z_p(z)$ and $Z_p(y)$ are transposed by the Hopf pairing, we have $(y_{m,l})= (R^{(c)}_{l,m}z_{l,m})$. By the proposition in section 1 we also know that \[Z_p(t_{y}(z))= \lbrace 2 \rbrace^{-1}(AZ_p(y)Z_p(z)- A^{-1}Z_p(z)Z_p(y)) \] so for all $m, l$ :  \[\begin{aligned}
z^{'}_{m,l}&=\lbrace 2 \rbrace^{-1}  \displaystyle { \sum_{j=0}^{N-1}} (Ay_{m,j}z_{j,l} -A^{-1}z_{m,j}y_{j,l})\\
&= \lbrace 2 \rbrace^{-1}  \displaystyle { \sum_{j=0}^{N-1}} (Az_{j,m} R^{(c)}_{j,m} z_{j,l} -A^{-1}z_{m,j} R^{(c)}_{l,j} z_{l,j})
\end{aligned}\]We have easily (by just writing the definition of $\lbrace  Q_n^{'(c)} \rbrace$) that for all $m, l$ : 

\[ z_{m,l}=\delta_{l,m-1}  \lbrace m \rbrace + \delta_{l,m} \lambda_{c+m} \]
Let us now compute $M_{m,l}^{(n)}= \lbrace n+1 \rbrace^{-1}(z^{'}_{m,l}- \lambda_{c+n} \delta_{m,l})$. We see from the expression of $(z_{m,l})$ that if $l \geq m+2$ or $m-2\geq l$ we have $M_{m,l}^{(n)}=0$ . Then we have :

\subparagraph{When $l=m-1$.} Then : \[\begin{aligned}
M_{m,m-1}^{(n)}&=(\lbrace 2 \rbrace \lbrace n+1 \rbrace)^{-1}  \displaystyle { \sum_{j=0}^{N-1}} (Az_{j,m} R^{(c)}_{j,m} z_{j,m-1} -A^{-1}z_{m,j} R^{(c)}_{m-1,j} z_{m-1,j})\\
&=(\lbrace 2 \rbrace \lbrace n+1 \rbrace)^{-1}  (Az_{m,m} R^{(c)}_{m,m} z_{m,m-1} -A^{-1}z_{m,m-1} R^{(c)}_{m-1,m-1} z_{m-1,m-1})
\end{aligned}\]
 Since $ R^{(c)}_{m,m}= R^{(c)}_{m-1,m-1}=1$ and  $ A \lambda_{c+m} - A^{-1} \lambda_{c+m-1}= (-A)^{-2N+2m}\lbrace 2 \rbrace $  (by computing), we have 
\[\begin{aligned}
M_{m,m-1}^{(n)}&= (\lbrace 2 \rbrace \lbrace n+1 \rbrace)^{-1} z_{m,m-1}( A \lambda_{c+m} - A^{-1} \lambda_{c+m-1})\\
& = \lbrace m \rbrace (-A)^{-2N+2m} (\lbrace n+1 \rbrace)^{-1}
\end{aligned}\]Then $M_{m,m-1}^{(n)}$ is clearly the evaluation at $X=A$ of a rational function $\widehat{M}_{m,m-1}^{(n)}(X) \in \mathbb{Q}(X)$ (which is clearly independent of $p$). And as $A \rightarrow -1$ : 
$$ M_{m,m-1}^{(n)} \rightarrow \dfrac{m}{n+1} $$
So $\widehat{M}_{m,m-1}^{(n)}(-1)=\dfrac{m}{n+1}$.

\subparagraph{When $l=m+1$.} A similar computation gives \[\begin{aligned}
M_{m,m+1}^{(n)}&=(\lbrace 2 \rbrace \lbrace n+1 \rbrace)^{-1} (Az_{m+1,m} R^{(c)}_{m+1,m} z_{m+1,m+1} -A^{-1}z_{m,m} R^{(c)}_{m+1,m} z_{m+1,m}) \\
&= (-A)^{-2N+2m+2} \lbrace -2N+2m+2 \rbrace \lbrace -2N+m+1 \rbrace^{+} \lbrace n+1 \rbrace^{-1}
\end{aligned}\]
Which tends to $\dfrac{-4(N-m-1)}{(n+1)}$ as $A \rightarrow -1$ . We see that $M_{m,m+1}^{(n)}$ can also be easily viewed as the evaluation at $X=A$ of a rational function $\widehat{M}_{m,m+1}^{(n)}(X) \in \mathbb{Q}(X) $ independent of $p$ which satisfies $\widehat{M}_{m,m+1}^{(n)}(-1)=\dfrac{-4(N-m-1)}{(n+1)}$.

\subparagraph{Finally when $l=m$.} 

\[\begin{aligned}
z^{'}_{m,m}&=\dfrac{Az_{m,m}^{2}R_{m,m}^{(c)}+Az_{m+1,m}^{2}R_{m+1,m}^{(c)}-A^{-1}z_{m,m-1}^{2}R_{m,m-1}^{(c)}-A^{-1}z_{m,m}^{2}R_{m,m}^{(c)}}{\lbrace 2 \rbrace}\\
&=\lbrace 2 \rbrace^{-1}((A-A^{-1}) \lambda_{c+m}^{2}+A \lbrace m+1 \rbrace^{2} R_{m+1,m}^{(c)} -A^{-1} \lbrace m \rbrace^{2} R_{m,m-1}^{(c)})
\end{aligned}\]
So  
$$M_{m,m}^{(n)}= \dfrac {(A-A^{-1}) \lambda_{c+m}^{2}+A \lbrace m+1 \rbrace^{2} R_{m+1,m}^{(c)} -A^{-1} \lbrace m \rbrace^{2} R_{m,m-1}^{(c)}}{\lbrace 2 \rbrace \lbrace  n+1 \rbrace}-\dfrac{\lambda_{c+n}}{ \lbrace n+1 \rbrace}$$
 Since for all $k$ , $\lambda_{c+k}$ is the evaluation at $X=A$ of $-((-X)^{-2N+2k+1}+(-X)^{2N-2k-1})$ (because $(-A)^ {p} = 1$ and $c=d-N$) which is independent of $p$ and by lemma 1 we deduce that there exists $\widehat{M}_{m,m}^{(n)}(X) \in \mathbb{Q}(X)$ (whose expression is clear)  independent of $p$ such that $\widehat{M}_{m,m}^{(n)}(A)= M_{m,m}^{(n)}$. An explicit check gives that as $A \rightarrow -1$ 
$$ \dfrac{ (A-A^{-1}) \lambda_{c+m}^{2} - \lambda_{c+n} \lbrace 2 \rbrace}  {\lbrace 2 \rbrace \lbrace  n+1 \rbrace} \rightarrow 0$$ 
And by lemma 1 as $A \rightarrow -1$ 
$$\dfrac{A \lbrace m+1 \rbrace^{2} R_{m+1,m}^{(c)}-A^{-1}\lbrace m \rbrace^{2} R_{m,m-1}^{(c)}}{ \lbrace 2 \rbrace \lbrace n+1 \rbrace} \rightarrow \dfrac{2(N-2m-1)}{n+1}$$
We can  conclude that $\widehat{M}_{m,m}^{(n)}(-1)= \dfrac{2(N-2m-1)}{n+1}$.   $\square$

\paragraph{Key observation : } The idea to prove theorem 1 is  very simple. Observe that if $n \leq N-2$ :
\[ Q_{n+1}^{'(c)}= \lbrace n+1 \rbrace^{-1}(Z_p(z)-\lambda_{c+n})( Q_{n}^{'(c)}) \]
 Since (see the proposition in section 1) \[ \rho_p^{(c)}(t_y)\bigl(Z_p(z)-\lambda_{c+n}Id \bigr) \rho_p^{(c)}(t_y)^{-1}= Z_p(t_y(z))-\lambda_{c+n}Id \]
 we have by inserting  $\rho_p^{(c)}(t_y)^{-1} \rho_p^{(c)}(t_y)$
\[\begin{aligned}
\rho_p^{(c)}(t_y)(Q_{n+1}^{'(c)}) &=\rho_p^{(c)}(t_y) \lbrace n+1 \rbrace^{-1}(Z_p(z)-\lambda_{c+n})( Q_{n}^{'(c)})\\
&=\rho_p^{(c)}(t_y) \lbrace n+1 \rbrace^{-1}(Z_p(z)-\lambda_{c+n}Id) \rho_p^{(c)}(t_y)^{-1} \rho_p^{(c)}(t_y)( Q_{n}^{'(c)}) \\
&=\lbrace n+1 \rbrace^{-1}(Z_p(t_y(z))-\lambda_{c+n})(\rho_p^{(c)}(t_y)(Q_{n}^{'(c)}))
\end{aligned}\]
In the basis $\lbrace Q_{n}^{'(c)} \rbrace$ this simply means that if we apply the matrix $M^{(n)}$ to the n-th column of the matrix $(a_{m,k})$ we get the (n+1)-st column of $(a_{m,k})$. In other words if we denote :
\[(a_{m,k})= (a_{0},...,a_{N-1}) \] where $a_i$ is the i-th column, we have
\[ a_{n+1}=M^{(n)}a_{n} \]
(Recall that $(a_{k,l})$ is the matrix of $\mu_c^{-1}\rho_p^{(c)}(t_y) $ in the basis $\lbrace Q_{n}^{'(c)} \rbrace$)

 From this key observation we are going to prove theorem 1 in 3 steps. First we  prove the existence of $\widehat{T}(X)$ , $\widehat{T}^{*}(X) \in GL_{N}(\mathbb{Q}(X)) $ independent of $p$ such that $\widehat{T}(A)=T_p$ and $\widehat{T}^{*}(A)=T_p^{*}$ ; then we  compute $\widehat{T}(-1)$ and $\widehat{T}^{*}(-1)$ ; finally we give an interpretation of $\widehat{T}(-1)$ and $\widehat{T}^{*}(-1)$.

\paragraph{Step 1 : Existence of $\widehat{T}(X)$  and $\widehat{T}^{*}(X)$ .} We define :
\[ e:=   \left( \begin{array}{c}
1 \\
0 \\
\vdots \\
0
\end{array} \right)  \in \mathbb{Q}^{N} \]
Let:
\[ \widehat{a}_{0}(X)=e\in \mathbb{Q}(X)^{N} \]
For $1 \leq n \leq N-1$, let:
\[   \widehat{a}_{n}(X)= \widehat{M}^{(n-1)}(X)...\widehat{M}^{(1)}(X) \widehat{M}^{(0)}(X)e \in \mathbb{Q}(X)^{N} \]
By lemma 2, these vectors are independent of $p$.
Since
\[ \rho_p^{(c)}(t_y)(Q_{0}^{'(c)})=\rho_p^{(c)}(t_y)(L_{c,0})=\mu_{c} L_{c,0}=\mu_{c}Q_{0}^{'(c)}  \]
we have
\[a_{0}= e =  \widehat{a}_{0}(A)\]
Then, by the key observation and lemma 2, for $1 \leq n \leq N-1$ :
\[\begin{aligned}
a_{n}&=M^{(n-1)}a_{n-1}\\
&=M^{(n-1)}...M^{(0)}a_{0} \\
&=\widehat{M}^{(n-1)}(A)...\widehat{M}^{(0)}(A)\widehat{a}_{0}(A) \\
&=\widehat{a}_{n}(A)
\end{aligned}\]
So $\widehat{T}(X):=(\widehat{a}_{m,n}(X))=(\widehat{a}_0(X),...,\widehat{a}_{N-1}(X)) \in GL_N(\mathbb{Q}(X)) $ is independent of $p$ and we have $(\widehat{a}_{m,n}(A))=(a_{m,n})$. On the other hand, 
\[\begin{aligned}
T^{*}_p&= (b_{n,m})\\
&= (a_{m,n} (R_{n,m}^{(c)})^{-1}) \\
&=( \widehat{a}_{m,n}(A)(\widehat{R}_{n,m}^{(c)}(A))^{-1}) \\
\end{aligned} \]
And by lemma 1, $\widehat{T}^{*}(X):= (\widehat{a}_{m,n}(X)(\widehat{R}_{n,m}^{(c)}(X))^{-1})  $ is independent of $p$. We have therefore found two matrices $\widehat{T}(X), \widehat{T}^{*}(X) \in GL_N(\mathbb{Q}(X))$ independent of $p$ such that $\widehat{T}(A)=T_p$ and $\widehat{T}^{*}(A)=T^{*}_p$.

\paragraph{Step 2 : Expression of $\widehat{T}(-1)$ and $\widehat{T}^{*}(-1)$.} We will prove that for all $n$:  
$$\widehat{a}_{m,n}(-1)=\dfrac{2^{n-m}(N-1-m)!}{(n-m)! (N-1-n)!} \ \text{when} \
m \leq n  \ \text{and}  \ 0 \ \text{otherwise.}\leqno(2)$$

 $$\widehat{b}_{n,m}(-1)=\dfrac{(-2)^{m-n}n!}{m!(n-m)!} \ \text{when} \ m\leq n \ \text{and} \ 0 \ \text{otherwise.}\leqno(3)$$
It is enough to prove $(2)$, since $(2)$ implies $(3)$ by using $\widehat{b}_{n,m}(-1) = \dfrac{\widehat{a}_{m,n}(-1)}{\widehat{R}_{n,m}(-1)}$ and lemma 1. We will compute $\widehat{T}(-1)$ by an induction on $n$ (the index of column).
 
If $n=0$,  $a_{0} =e $ so the limit is as expected.

If $n=1$, by the key observation \[ \widehat{a}_{1}=\widehat{M}^{(0)}\widehat{a}_{0}\]
so by lemma 2 : 
\[ \begin{aligned}
\widehat{a}_{0,1}(-1)&= \widehat{M}^{(0)}_{0,0}(-1)=2(N-1) \\
\widehat{a}_{1,1}(-1)&= \widehat{M}^{(0)}_{1,0}(-1)=1 \\
\widehat{a}_{m,1}(-1)&=0 \  \ \ \text{when}  \ m >1
\end{aligned} \]
So $(2)$ is true when $N=2$. 

Now suppose $N \geq 3$. Let  $1 \leq n \leq N-2 $ and suppose by induction that $(2)$ holds for $n$.
Then by the key observation$$\widehat{a}_{n+1}(-1)= \widehat{M}^{(n)}(-1)\widehat{a}_{n}(-1)$$ 
So by lemma 2, when $1 \leq m \leq n-1$ :
 \[\begin{aligned}
\widehat{a}_{m,n+1}(-1) &=  \widehat{M}^{(n)}_{m,m-1}(-1) \widehat{a}_{m-1,n}(-1) +\widehat{M}^{(n)}_{m,m}(-1) \widehat{a}_{m,n}(-1) +  \widehat{M}^{(n)}_{m,m+1}(-1) \widehat{a}_{m+1,n}(-1) \\
&= \dfrac{m}{n+1} \dfrac {2^{n+1-m}(N-m)!}{(n+1-m)! (N-n-1)!} +\dfrac{2(N-2m-1)}{n+1} \dfrac {2^{n-m}(N-m-1)!}{(n-m)! (N-n-1)!}  \\
& \ \ \ \  \ \ \ \ \ \ \ \ \  \ \ \ \ \ \ \ \ \  \ \ \ \ \ \ \ \ \  \ \ \ \ \ \ \ \ \  \ \ \ \ \   \ -\dfrac{4(N-m-1)}{n+1} \dfrac {2^{n-1-m}(N-m-2)!}{(n-1-m)! (N-n-1)!} \\
&  =\dfrac {2^{n+1-m}(N-m-1)!}{(n+1-m)! (N-n-2)!}\Biggl( \dfrac{m(N-m)}{(n+1)(N-n-1)} \\
&\ \ \ \  \ \ \ \ \ \ \ \ \  \ \ \ \ \ \ \ \ \  \ \ \ \ \ \ \ \ \  \ \ \ \ \ \ \ \ \ \ +\dfrac{(N-2m-1)(n+1-m)}{(n+1)(N-n-1)}-\dfrac{(n-m)(n+1-m)}{(n+1)(N-n-1)}\Biggr) \\
&= \dfrac {2^{n+1-m}(N-m-1)!}{(n+1-m)! (N-n-2)!}
\end{aligned}\]
It remains to consider the cases $m=0$, $m=n$, $m=n+1$	and $m \geq n+2$. We leave it to the reader to check that:
\[ \begin{aligned}
& \widehat{a}_{0,n+1}(-1) = \dfrac{2^{n+1} (N-1)!}{(n+1)! (N-n-2)!} \\
& \widehat{a}_{n,n+1}(-1) =2(N-n-1) \\
& \widehat{a}_{n+1,n+1}(-1) = 1 \\
& \widehat{a}_{m,n+1}(-1)=0 \ \text{when} \ m \geq n+2
\end{aligned} \]So we have shown that $(1)$ holds for $n+1$. This completes the proof of the induction step.

\paragraph{Step 3 : Interpretation of $\widehat{T}(-1)$ and $\widehat{T}^{*}(-1)$ .} Recall the action of $SL_{2}( \mathbb{Z})$ on $H_N$ (the space of homogeneous polynomials of two variables X and Y of total degree $N-1$) . 
For $0 \leq n \leq N-1$ and $ \begin{pmatrix}
a& b \\
c&d
\end{pmatrix} \in SL_{2}( \mathbb{Z})$ :
 
\[ \begin{pmatrix}
a& b \\
c&d
\end{pmatrix}(X^{N-n-1}Y^{n}) = (aX+cY)^{N-n-1}(bX+dY)^{n} \]
so 

$ \begin{pmatrix}
1& 1 \\
0&1
\end{pmatrix}(X^{N-n-1}Y^{n}) =  \displaystyle { \sum_{m=0}^{n}} \dfrac{n!}{m!(n-m)!} X^{N-m-1}Y^{m} $

$\begin{pmatrix}
1& 0 \\
-1&1
\end{pmatrix}(X^{N-n-1}Y^{n}) =  \displaystyle { \sum_{m=n}^{N-1}} \dfrac{(-1)^{m-n} (N-n-1)!}{(m-n)!(N-m-1)!} X^{N-m-1}Y^{m} $

It gives in the basis  $(1)$ (see theorem 1) :
\begin{itemize}
\item $ \begin{pmatrix}
1& 1 \\
0&1
\end{pmatrix}(\alpha_{n} X^{N-n-1}Y^{n}) =\displaystyle { \sum_{m=0}^{N-1}}\widehat{a}_{m,n}(-1) (\alpha_{m} X^{N-m-1}Y^{m}) $

\item $\begin{pmatrix}
1& 0 \\
-1&1
\end{pmatrix}(\alpha_{n} X^{N-n-1}Y^{n}) =\displaystyle { \sum_{m=0}^{N-1}}\widehat{b}_{m,n}(-1) (\alpha_{m} X^{N-m-1}Y^{m}) $

\end{itemize}

Since $t_{y}$ maps to $ \begin{pmatrix}
1& 1 \\
0&1
\end{pmatrix}$ and $t_{z}$ maps to $\begin{pmatrix}
1& 0 \\
-1&1
\end{pmatrix}$ in $SL_2(\mathbb{Z})$ we conclude that $\widehat{T}(-1)$ and $\widehat{T}^{*}(-1)$ are the matrices of $h_{N}(t_{y})$ and $h_{N}(t_{z})$ in the basis $(1)$ which completes the proof of theorem 1. $\square$

\paragraph{Remark 1 :}Using the previous techniques, one can get explicit formulas for $\widehat{T}$ and   $\widehat{T}^{*}$  but they are quite complicated and we don't need them to compute the limits.

\paragraph{Remark 2 :} The main theorem is an easy implication of theorem 1. Indeed for all odd $p$  and $ p \geq 2N+1$ , $ \rho_{p}^{(c_p)}$ is a representation so :

\[ \widehat{T}(A) \widehat{T}^{*}(A)  \widehat{T}(A)  = \widehat{T}^{*}(A) \widehat{T}(A)  \widehat{T}^{*}(A) \](for any  primitive 2p-th root of unity $A $). Since a rational function has only a finite number of roots, in $GL_{N}(\mathbb{Q}(X))$ we have :

\[ \widehat{T} \widehat{T}^{*}  \widehat{T}  = \widehat{T}^{*} \widehat{T} \widehat{T}^{*}  \]
Since $\Gamma_{1,1}$ has a presentation $\langle t_y, t_z \mid t_y t_z t_y = t_z t_y t_z \rangle $, the previous relation ensures that there exists a unique representation $ \rho : \Gamma_{1,1} \rightarrow GL_{N}(\mathbb{Q}(X))$ such that $\rho(t_{y})=\widehat{T}$ and $\rho(t_{z})=\widehat{T}^{*}$. By the same argument, for all odd $p$ and all primitive 2p-th roots of unity $A$, there exists a unique character $\chi_p : \Gamma_{1,1} \rightarrow \mathbb{C}^{*}$ such that $\chi_p(t_{y})= \mu_{c_p}$ and $\chi_p(t_{z})= \mu_{c_p}$. Then if we choose the same $A$ (a primitive 2p-th root of unity) to define  $\rho_p^{(c_p)}$ and $\chi_p$, by theorem 1 in the bases considered above we have:

\[ \chi_p  \otimes \rho^{[A]} = \rho_p^{(c_p)}\]

and 

\[  \rho^{[-1]}= h_N \]

\subsection{Proof of theorem 2}

We know that $\Gamma_{1,1}$ acts on $H_1(\mathbb{T}^{2}, \mathbb{C})$ (the first homology of the torus with coefficients in $\mathbb{C}$). So we have a representation : \[ \varphi_{N} : \Gamma_{1,1} \longrightarrow Aut(Sym^{N-1}H_1(\mathbb{T}^{2}, \mathbb{C})) \]where $Sym^{N-1}H_1(\mathbb{T}^{2}, \mathbb{C})$ is the space of  symmetric (N-1)-tensors on $H_1(\mathbb{T}^{2}, \mathbb{C})$. We then use the following fact: $\varphi_N$ is isomorphic to $h_N$ (because $\varphi_2$ is isomorphic to $h_2$).

Now let $\phi \in \Gamma_{1,1}$ be a pseudo-Anosov, we denote by $\phi'$ the image of $\phi$ in $SL_{2}(\mathbb{Z})$. We have $\mid tr( \phi') \mid  > 2$  so there exists $u  \in H_1(\mathbb{T}^{2}, \mathbb{C}) - \lbrace 0 \rbrace$ and $\lambda \in \mathbb{C}$ with $\mid \lambda \mid > 1$ such that $\phi' u = \lambda u $ . Then $u^{\otimes N-1} \in Sym^{N-1}H_1(\mathbb{T}^{2}, \mathbb{C})$ and  $h_{N}(\phi) u^{\otimes N-1} = \lambda^{N-1} u^{\otimes N-1}$.

Let $M$ be the matrix of $h_N(\phi)$ in the basis  of theorem 1 . We deduce that $M$ has an eigenvalue $\lambda^{N-1}$ such that $\mid \lambda^{N-1} \mid > 1 $.

Now for odd $p$  and $p \geq 2N+1$ we set $A_p= - e^{\frac{i \pi}{p}}$  which is a primitive 2p-th root of unity. We define $\rho_{p}^{(c_p)}$ using this root $A_p$. Let $M_p$ be the matrix of $\rho_{p}^{(c_p)}(\phi)$ in the basis $ \lbrace Q_n^{'(c)} \rbrace$. By the main theorem since $ A_p \rightarrow -1 $ (as $p \rightarrow  \infty$):
\[ (\chi_p(\phi))^{-1} M_{p} \rightarrow M \ \ \text{when} \ p \rightarrow  \infty \]So there exists $p_0(\phi)$ such that for $p \geq p_0(\phi)$ , the matrix $(\chi_p(\phi))^{-1} M_{p}$ has an eigenvalue whose absolute value is greater than 1 and hence this matrix has infinite order. Since $\chi_p(\phi) $ is a root of unity for all $p$, it follows that the matrix $M_p$ also has infinite order for $ p \geq p_0(\phi)$. Since $M_p$ is the matrix of $\rho_{p}^{(c_p)}(\phi)$, theorem 2 is proved.

\paragraph{Remark : } This proof is very similar to the case of the sphere with 4 holes see \cite{1}.

 Institut de Math\'ematiques de Jussieu (UMR 7586 du CNRS)

Equipe Topologie et G\'eom\'etrie Alg\'ebriques,

Case 247, 4 pl.Jussieu,

75252 Paris Cedex 5, France.

Email: santharo@clipper.ens.fr


\begin{thebibliography}{2}





\bibitem[1]{1} Joergen E. Andersen, Gregor Masbaum, Kenji Ueno. Topological Quantum Field Theory and the Nielsen-Thurston classification of M(0,4). Math. Proc. Cam. Phil. Soc. 141 (2006) 477-488 

\bibitem[2]{2} C. Blanchet, N. Habegger, G. Masbaum, P. Vogel. Topological quantum field theories derived from the Kauffman bracket, Topology 34 (1995) 883-927

\bibitem[3]{3}P. M. Gilmer, G. Masbaum. Integral TQFT for a one-holed torus. Pac. J. Math. 252 (2011) No. 1, 93-112

\bibitem[4]{4} G. Masbaum, J. Roberts. On central extensions of mapping class groups, Math.Ann 302 (1995), 131-150. MR 96i: 57013













\end{thebibliography}
\end{document}